\documentclass[12pt,righttagstyle]{amsart}

\usepackage{amssymb,latexsym}

\def\cal{\mathcal}
\def\Bbb{\mathbb}

\newcommand{\bC}{{\Bbb C}}

\newcommand{\bR}{{\Bbb R}}

\def\pf{\proof}
\def\bpf{\pf}
\newcommand{\epf}{ {$\Box$}\medskip}
\def\a{\alpha}
\def\b{\beta}

\def\Ga{\Gamma}
\def\c{\xi}
\def\d{\delta}

\def\w {\omega}
\def\L{\Lambda}
\def\s{\sigma}
\def\x{\chi}
\def\p{\pi}

\def\n{\nu}
\def\r{\rho}
\def\t{\tau}
\def\f{\phi}
\def\F{\Phi}
\def\th{\theta}

\def\et{\eta}
\def\z{\zeta}
\def\W{\Omega}
\def\De{\Delta}
\def\del{\nabla}

\def\C{\bold C}
\def\HH{\text{\bf H}}

\def\R{\bold R}

\def\Ad{\text{Ad}}

\def\det{\text{det}}
\def\Tr{\text{Tr}}
\def\sx{\cal X}
\def\sy{\cal Y}

\def\ov{\overline}
\def\p{\partial}

\def\sd{\cal D}

\def\sh{\cal H}

\def\sj{\cal J}

\def\sL{\cal L}

\def\sn{\cal N}

\def\sp{\cal P}
\def\sq{\cal Q}

\def\sx{\cal X}
\def\sy{\cal Y}
\def\sz{\cal Z}
\def\dP{\Psi _*}

\def\w0{w^{(0)}}

\def\C{{\cal C}}

\def\qed{{\scriptstyle{\Box}}}
 
\def\De{\Delta}
\def\del{\nabla}

\def\C{\bold C}
\def\HH{\text{\bf H}}

\def\R{\bold R}

\def\Ad{\text{Ad}}

\def\det{\text{det}}
\def\Tr{\text{Tr}}
\def\nfi{N(\Phi)}
\def\<{\langle}
\def\>{\rangle}
\newtheorem{thm}{Theorem}[section]
\newtheorem{prop}[thm]{Proposition}
\newtheorem{cor}[thm]{Corollary}
\newtheorem{lem}[thm]{Lemma}

\begin{document}

\title[ Hua system and pluriharmonicity]{ Hua system and
pluriharmonicity
for symmetric irreducible Siegel domains of type II }

\author[Bonami]{Aline Bonami$^1$}
\address{MAPMO-UMR 6628,
D\'epartement de Math\'ematiques, Universit\'e d'Orl\'eans, 45067
Orl\'eans Cedex 2, France} \email{{\tt
Aline.Bonami@labomath.univ-orleans.fr}}
\author[Buraczewski]{Dariusz Buraczewski$^{1,2}$}
\address {Instytut Matematyczny, Uniwersytet Wroc\l awski,
Plac Grunwaldzki 2/4, 50-384 Wroc\l aw, Poland} \email
{dbura@math.uni.wroc.pl }
\author[Damek]{Ewa Damek$^{1,2}$}
\address {same address in Wroc\l aw}
\email {{\tt edamek@math.uni.wroc.pl}}
\author[Hulanicki]{Andrzej Hulanicki$^{1,2}$}
\address {same address in Wroc\l aw}
\email {{\tt hulanick@math.uni.wroc.pl}}
\author[Penney]{Richard Penney}
\address {Department of Mathematics, Purdue University, West
Lafayette, Indiana 47907, USA} \email {{\tt rcp@math.purdue.edu} }
\author[Trojan]{Bartosz Trojan$^{1,2}$}
\address {same address in Wroc\l aw}
\email {{\tt s85183@math.uni.wroc.pl} }
\thanks{ $^1$ This work was  partly done within the project TMR
Network
``Harmonic Analysis",  contract no. ERB
FMRX-CT97-0159.
We thank the European Commission and the mentioned Network
for the support provided.\\
$^2$ The authors were
partly supported by KBN grant 2P03A04316, Foundation for Polish
Sciences,
Subsidy 3/99}

\begin{abstract}
    We consider here a generalization of the Hua system which was
    proved by Johnson and Kor\'anyi to characterize Poisson-Szeg\"o integrals for Siegel
    domains of tube type. We show that the situation is
    completely different when dealing with non tube type symmetric
    irreducible symmetric domains: then all functions which are
    annihilated by this second order system and satisfy an $H^2$ type
    integrability condition are pluriharmonic functions.

\end{abstract}
\maketitle

\section{Introduction}

Let $\sd$ be a bounded symmetric domain in $\C^m$, and let $G$ be
the group of all biholomorphic transformations of $\sd$. The aim of
this paper is to study $\HH$-harmonic functions, where $\HH$ is a
naturally
defined $G$-invariant real system of second order
differential
operators on $\sd$ which annihilates  pluriharmonic functions.
The system $\HH$ is defined in terms of the K\"ahler structure of
$\sd$ and makes sense on every K\"ahlerian manifold.

\medskip
To define the system $\HH$, we recall some basic facts about $\sd
$.
Let $T^{1,0}(\sd)$ be the holomorphic tangent bundle of $\sd$. The
Riemannian connection $\del $ induced by the Bergman metric on $\sd
$
preserves $T^{1,0}(\sd)$ and so does the curvature tensor.
For $Z,W$ two complex vector fields we denote by
  $R(Z,W)=\del_Z\del_W-\del_W\del_Z-\del_{[Z,W]}$
the  curvature tensor restricted to $T^{1,0}(\sd)$.
Let $f$ be a smooth function on $\sd$ and let
\begin{equation}
  \De(Z,W)f=(Z\overline W-\del_{Z}\overline W)f= (\overline W
Z-\del_{\overline
  W}Z)f\,.
    \label{delta}
\end{equation}
Then $\De(Z,W)$ may be seen as a second order operator which
annihilates both holomorphic and antiholomorphic functions, and
consequently, the pluriharmonic functions.
Conversely, if all $\De(Z,W)$ annihilate $f$, then $f$ is
pluriharmonic. Indeed, we have
$\De(\partial_{z_{j}},\partial_{z_{k}})=
\partial_{z_{j}}\partial_{\overline{z_{k}}}$.

Let $(\ ,\ )$ be the canonical Hermitian product in
$T^{1,0}(\sd )$.
Fixing a smooth function $f$, we use $(\ ,\ )$ to define a smooth
section $\Delta _f$ of the bundle of endomorphisms of
$T^{1,0}(\sd)$:
\begin{equation}
    (\Delta f\cdot Z, W)=\Delta (W,Z)f\, ,
    \end{equation}
where $Z,W$ are holomorphic vector fields. Then
we define $\HH f$ as another smooth section of the bundle of
endomorphisms of
$T^{1,0}(\sd )$ by
\begin{equation}
    (\HH f\cdot Z, W)=\Tr(R(\overline
Z,W)^*\De_{f})=\Tr(R(\overline W,Z)
\De_{f})\,.
    \label{Hua-inv}
\end{equation}
To compute explicitly $\HH f$, we may take an orthonormal frame of
sections of $T^{1,0}(\sd)$, which we denote $E_1,E_2,\cdots,E_m$.
Then
\begin{equation}
    \HH f
    =\sum _{j,k}(\De (E_j, E_k)f) R(\overline E_{j}, E_{k})
    \label{Hua-inv2}\,.
\end{equation}
The system $\HH$ is, of course, a contraction of the tensor
field $\De_{f}$. It is invariant with respect to biholomorphisms,
 which means that
\begin{equation}
\HH (f\circ \Psi)=\dP ^{-1}[(\HH f)\circ \Psi]\dP
\label{inv}
\end{equation}
for every biholomorphic transformation $\Psi $ of $\sd$, $\dP$
being its  differential.

By definition, $\HH$-harmonic functions are functions  which are
annihilated by $\HH$. We will consider here symmetric Siegel
domains, for which these notions are well defined since they are
holomorphically equivalent to bounded domains. When $\sd$ is a
symmetric Siegel domain of {\it tube type}, (\ref{Hua-inv}) is
equivalent to the classical Hua system.  This system is known to
characterize the Poisson-Szeg\"o integrals (see \cite {FK} and
\cite {JK}). This means that a function on $\sd$ is $\HH$-harmonic
if, and only if, it is the Poisson-Szeg\"o integral of a
hyperfunction on the Shilov boundary. Originally, the curvature
tensor was not explicit in the Hua system. For  classical domains,
the system has been defined by L. K.  Hua as a ``quantization'' of
the equation defining the Shilov boundary (see \cite {Hua} and
\cite {BV}). L. K. Hua proved that the system annihilates
Poisson-Szeg\"o integrals. Then the system was extended by K.
Johnson and A. Kor\'anyi \cite {JK} to all symmetric tube type
domains and was written down in terms of the enveloping algebra of
the semi-simple Lie group of automorphisms of the domain. K.
Johnson and A. Kor\'anyi  proved  not only that for all tube
domains the system annihilates the Poissonn-Szeg\"o kernel, but
also that the $\HH$-harmonic functions are Poisson-Szeg\"o
integrals. Rewriting Johnson--Kor\'anyi formula $C(\partial ,\bar
\partial )$ in terms of the curvature tensor,  as suggested by
Nolan Wallach, one obtains the same system as above. It is why we
call $\HH$ the Hua-Wallach system.

Notice that (\ref{Hua-inv2}) and (\ref{inv})  have a perfect sense
on any K\"ahlerian manifold, and, for general Siegel domains, the
system (\ref{Hua-inv2}) has been already studied in \cite {DHP}.
In particular, for non tube symmetric Siegel both (\ref{Hua-inv2})
and Johnson--Koranyi formula $C(\partial ,\bar
\partial )$
take the same form. In the work of N. Berline and M. Vergne \cite
{BV} it is observed that $C(\partial ,\bar \partial )$ does not annihilate
Poisson - Szeg\"o integrals, and  the problem of
describing
$C(\partial ,\bar \partial )$--harmonic functions is risen. Here we are going
to answer
their question.

\medskip

\noindent {\bf Main Theorem}. {\it Let $\sd$ be a symmetric
irreducible Siegel domain of type II, and let $F$ be an
$\HH$-harmonic function on $\sd$ which satisfies the growth
condition}
\begin{equation}\label{H2}
\tag{$H^2$}
\sup_{z\in\sd}\int_{N(\Phi)}\vert F(uz)\vert^2 du<\infty ,
\end{equation}
{\it where $N(\F )$ is a nilpotent subgroup of $S$ whose action is
parallel to the Shilov boundary. Then $F$ is pluriharmonic.}

\medskip

This is in a striking contrast to the case when $\sd $ is a
symmetric
tube domain. It requires some comments.

The Poisson-Szeg\"o
integrals
on type II domains have been characterized by N. Berline and
M. Vergne \cite{BV} as zeros of a $G$-invariant system
which, in general, is of the third order. It is obtained by
``quantization''
of the Shilov boundary equations. They also prove that for
 domains over the cone of
hermitian
positive definite matrices   one can use a second order system, $\De _Z$,
to characterize Poisson-Szeg\"o integrals. This  system
 appears already in the book by Hua \cite {Hua}. It is
obtained from
$C(\partial ,\bar \partial )$ by a projection that eliminates a part of the
equations.

All this shows that the system  $\HH$  does not seem
to be canonical in any sense, although it is
defined with the aid of the curvature tensor, certainly an
important invariant, the geometric
meaning of the system being
still unclear. Our present work suggests that it would be
interesting
to understand second order  systems of operators on
symmetric Siegel domains which are invariant under the full group of
biholomorphisms.

In the proof of the main theorem, we use heavily the theory of
harmonic functions
with respect to subelliptic operators
on solvable Lie groups  \cite {R},  \cite {D},  \cite {DH},  \cite
{DHP}.
To do this, we identify the
domain
$\sd$ with a solvable Lie group $S\subset G$ that acts simply
transitively
on $\sd$. We then use a special orthonormal frame  of $S$-invariant
vector
fields, $E_1,E_2,\cdots,E_m$, to compute the operator $\HH$ by the
formula
(\ref{Hua-inv2}). In fact, we only consider the left-invariant
second order elliptic
operators built
out of the diagonal of $\HH$,
\begin{equation}
    \HH _{j}f=(\HH f\cdot E_j, E_j).
    \label{coord}
\end{equation}
Elliptic operators which are linear combinations of
operators $ \HH _{j}$
play the main role in our argument, and in particular the
Laplace-Beltrami operator $\De$, which is the trace of $\HH$.
We represent $\HH$-harmonic
functions
as various Poisson integrals, and we use properties of these
representations.

Linear combinations of the operators $\De (E_j, E_{k})$ have
already
been used to characterize pluriharmonic functions (see
\cite{DHMP}).
We should emphasize that the systems under
study here are different from those of \cite{DHMP}, and the proofs
require new ideas. Since a part of the construction is the same in
the
two papers, we try to simplify the presentation for the
reader's convenience.

Our growth assumption $(H^2)$ is made mainly  for
technical reasons, $L^2$ harmonic analysis being the
easiest. We hope to be
able to obtain similar conclusions for bounded functions, and
perhaps even for larger classes of functions. This requires,
however, a somewhat more delicate technic.
In fact, it is not clear that the conclusion requires any growth
condition at the
boundary, and one may conjecture that only growth conditions at
infinity
are necessary to insure pluriharmonicity for $\HH$--harmonic
functions on symmetric irreducible Siegel domains of type
II. On the other hand, for tube type domains, one may conjecture
that
growth conditions on derivatives at the boundary insure
pluriharmonicity
for $\HH$-harmonic functions as it is the case in the unit ball
(\cite{BBG}) for $\De$-harmonic functions.

Finally, let us remark that, if we do not insist on invariance
properties of the systems considered, then it is always possible
to characterize pluriharmonic functions, among the functions which
are harmonic with respect to the Laplace-Beltrami operator
$\Delta$, as those which are annihilated by a single second order
operator $L$ (without any growth condition). Indeed,  a classical
theorem of Forelli (see Rudin's book \cite{Ru}) asserts that every
smooth function in the unit ball which is annihilated by the
operator $\sum z_{j} \overline {z_{k}}\frac{\partial^2} {\partial
z_{j}\partial \overline {z_{k}}}$ is pluriharmonic in the ball. So
$L$ can be taken as this operator suitably translated, so that a
function which is annihilated by $L$ is pluriharmonic in the
neihborhood of a point. Then the real-analyticity of the function,
which follows from the fact that it is $\Delta$--harmonic, insures
its pluriharmonicity everywhere.

In view of  Forelli's Theorem,  it is not so much the small number
of  operators in the system used to characterize pluriharmonic
functions than the strong invariance properties of the system
itself which are relevant. In this context,
 the present paper can be viewed as a complement to \cite{DHP}
and\cite{DHMP}.

\medskip

\section {Hua-Wallach systems}
  \subsection {General Hua-Wallach systems}
In this subsection, $\sd$ is a general domain in $\bC^m$ which is
holomorphically equivalent to a bounded domain. We recall here the
properties of the K\"ahlerian structure related to the Bergman
matric as well as
some elementary facts about the Hua-Wallach system which we will
use later.
The reader may refer to
\cite {H} and \cite {KN} for more details on the prerequisites.

Let $T$ be the tangent bundle for the complex domain $\sd$, and let
$T^{\C}$ be the
complexified
tangent bundle.
The complex structure $\sj $ and the Bergman metric $g$ are
extended from $T$ to
$T^{\C}$ by complex
linearity.
Let $T^{1,0}$ and $T^{0,1}$ be the eigenspaces of $\sj $ such that
$\sj
|_{T^{1,0}}=
i \text{Id}$, $\sj |_{T^{0,1}}=-i \text{Id}$. We have
$$T^{\C}=T^{1,0}\oplus T^{0,1}\,.$$
The conjugation operator exchanges $T^{1,0}$ and $T^{0,1}$.

The spaces of smooth
sections of $T$, $T^{\C}$, $T^{1,0}$
will be
denoted $\Ga (T)$, $\Ga (T^{\C})$, $\Ga (T^{1,0})$, respectively.
Smooth sections of $T^{1,0}$ are called holomorphic vector fields.

 The Riemannian connection $\del $ is also extended from $\Ga (T)$
to $\Ga(T^{\C})$ by
complex
linearity and, since $\del $ is defined by a K\"ahlerian structure,
it commutes
with $\sj$. An immediate consequence is that $\del _Z W$ belongs to
$\Ga(T^{1,0})$
(respectively $\Ga(T^{0,1})$) whenever $W\in \Ga(T^{1,0})$
(respectively
$\Ga(T^{0,1})$). Moreover, for every couple $U,V \in \Ga(T^{\C})$,
we
have that $\ov{\del _U V}=\del _{\ov U} \ov V$ and
\begin{equation}
[\overline U,V]= \del_{\overline U}V-\del _V {\overline
U},\label{brack}.
\end{equation}
Therefore, for $Z, W$ holomorphic vector fields,
\begin{equation}
   \del_{\overline W}Z = \pi_{(1,0)}([\overline W,Z]) ,
    \label{connec}
\end{equation}
where  $\pi_{(1,0)}$ denotes the projection from $T^{\C}$ onto
$T^{1,0}$.
 The curvature
tensor $$
R(U,V)=\del _U\del _V-\del _V\del _U-\del _{[U,V]}
$$
preserves also  $T^{1,0}$ and $R(\ov U,\ov V)\ov
Z=\ov{R(U,V)Z}$.
The restriction of $R(U,V)$ to $T^{1,0}$ is also denoted by
 $R(U,V)$. On $T^{1,0}$ the Hermitian
scalar product arising from the Bergman metric is denoted by
$$
(Z,W)=\frac{1}{2}g(Z, \ov W).
$$
For $U$, $V$, $Z$, $W$ holomorphic vector fields, we have
\begin{equation}
(R( \ov V,U)Z,W)=(R(\ov W,Z)\, U, V)=( U,R(\ov Z,W)\, V)\,.
\label{curv}
\end{equation}
In particular,
\begin{equation}
R(\ov W,Z)=R(\ov Z,W)^*\,.\label{star}\end{equation} Let us now go
back to the definitions given in the introduction. The identity $
Z\overline W-\del_{Z}\overline W= \overline W Z-\del_{\overline
  W}Z$ is a direct consequence of (\ref{brack}). The fact that all
$\De(Z,W)$
  annihilate  pluriharmonic functions follows from (\ref{connec})
as well as from the identity
$\De(\partial_{z_{j}},\partial_{z_{k}})=
\partial_{z_{j}}\partial_{\overline{z_{k}}}$. Moreover,
  $$\De(\phi Z, \psi W)f=\phi\,\ov\psi \,\De(Z,W)f\,.$$
which means that $\De_{f}$ is a tensor field.
The equality in (\ref{Hua-inv}) comes from (\ref{brack}),
while one proves (\ref{Hua-inv2})
 using (\ref{curv}) for $(R(\ov W,Z)\cdot E_k,E_j)$.

Let us now show invariance of $\HH_{\sd}$ with respect to
biholomorphisms.
Let $\Psi $ be a biholomorphism from $\sd$ onto $\sd'$,
and $\Psi _*$ the  holomorphic differential of $\Psi $ which
maps $T^{1,0}_\sd$ onto $T^{1,0}_{\sd'}$. All tensor fields are
transported by $\Psi$, including, of course,
the Riemannian structure and the curvature tensor.
Thus
 $$R_\sd (\ov W,Z)=\Psi_*^{-1}R_{\sd'} (\ov{
\Psi_*W},\Psi_*Z)\Psi_*\,.$$
Moreover, for a smooth function $f$ on $\sd'$ and $g=f\circ \Psi$,
we have
$\De_g= \Psi_*^{-1} \De_f\Psi_*$. So
$$(\HH_\sd \,g\cdot Z,W)=\Tr\left(R_{\sd'} (\ov{ \Psi_* W},\Psi_*Z)
\De_f\right)$$
 which implies
$$\HH_\sd \,g=\Psi_*^{-1} (\HH_{\sd'} \,f) \Psi_*\,.$$

Finally, let us remark that, from formulas (\ref{delta}),
(\ref{Hua-inv}) and (\ref{star}),
it follows that $\De_{\ov f}=(\De_{f})^*$, and $\HH{\ov f}=(\HH
f)^*$.
So, to study $\HH$-harmonic functions, it is sufficient to consider
functions
which are real-valued.

\medskip
We want now to compute explicitly the Hua-Wallach operator
 for symmetric irreducible Siegel domains. To do it, we will use
Formula (\ref{Hua-inv2})
for a particular orthonormal basis $E_1,...,E_m$.

\medskip

\subsection {Preliminaries on irreducible symmetric cones}

Let $\Omega$ be an irreducible symmetric cone in an Euclidean
space.
Our aim is to describe
precisely the solvable group  that acts simply
transitively on $\Omega$. The group will be used
 in the
construction
of the orthonormal basis. We do it all in terms of Jordan
algebras, and we refer to the book of
Faraut and Kor\'anyi \cite {FK}
for these prerequisites, introducing here only
 the notations and  principal results that  will be needed later.

A finite dimensional algebra $V$ with a scalar product
$\langle\cdot ,\cdot \rangle$ is an Euclidean Jordan algebra if
for all elements $x, y$ and $z$ in $V$ $$xy=yx \hspace{2cm}x(x^2y)
= x^2(xy)\hspace {2cm} \langle xy,z\rangle=\langle y,xz\rangle. $$
 We denote by $L(x)$ the
self-adjoint endomorphism of $V$ given by the multiplication by
$x$, i.e. $L(x)y=xy.$

For an irreducible symmetric cone $\Omega$ contained in a linear
space $V$ of same dimension, the space $V$ can be made a simple
real Euclidean Jordan algebra with unit element $e$, so that
$$\Omega=\,\mbox{int}\,\{x^2:\ x\in V\}.$$ Let $G$ be the
connected component of the group  of all transformations in
$GL(V)$ which leave $\Omega$ invariant, and let ${\cal G}$ be its
Lie algebra. Then ${\cal G}$ is a subspace of the space of
endomorphisms of $V$ which contains all $L(x)$ for all $x\in V$,
as well as all $x \qed y$ for $x,y \in V$, where $x \qed
y=L(xy)+[L(x),L(y)]$ (see \cite{FK} for  these properties).

\smallskip

We fix a Jordan frame $\{c_1,\dots,c_r\}$ in $V$, that is, a
complete system of orthogonal primitive idempotents: $$
c^2_i=c_i,\hspace {2cm} c_ic_j=0 \quad\mbox{if}\; i\neq j, \hspace
{2cm}c_1+...+c_r=e $$ and none of the $c_1,...,c_r$ is a sum of
two non-zero idempotents. Let us recall that the length $r$ is
independent of the choice of the Jordan frame. It is called the
rank of $V$. To have an example in mind, one may think of the
space $V$ of the symmetric $r\times r$  matrices endowed with the
symmetrized product of matrices $\frac{1}{2}(xy+yx)$. Then the
corresponding cone is the set of symmetric positive definite
$r\times r$  matrices, the set of diagonal matrices with all
entries equal to $0$ except for one equal to $1$ being a Jordan
frame.

The Peirce decomposition of $V$
related to the Jordan frame
$\{c_1,\dots,c_r\}$ (\cite {FK}, Theorem IV.2.1) may be written as
\begin{equation}
    V=\bigoplus_{1\leq i\leq j\leq r}V_{ij}\,.
    \label{Peirce}
\end{equation}
It is given by the common diagonalization of the self-adjoint
endomorphims $L(c_{j})$ with respect to their only eigenvalues
$0$, $\frac 12$, $1$. In particular $V_{jj}=\bR c_j$ is the
eigenspace of $L(c_{j})$ related to $1$, and, for $i<j$, $V_{ij}$
is  the intersection of the eigenspaces of $L(c_{i})$ and
$L(c_{j})$ related to $\frac 12$. All $V_{ij}$, for $i<j$, have
the same dimension $d$.

 For each $i<j$, we fix once for all an
orthonormal basis of $V_{ij}$, which we note
 $\{e^\alpha_{ij}\}$, with $1\leq \alpha\leq d$. To simplify the
notation,
 we write $e^\alpha_{ii}=c_{i}$ ($\alpha$ taking only the
 value $1$). Then the system  $\{e^\alpha_{ij}\}$, for $i\leq j$
 and $1\leq \alpha \leq \dim V_{ij}$,
 is an orthonormal basis of $V$.

Let us denote by $\cal A$ the  abelian subalgebra of $\cal G$
consisting of elements
$H=L(a)$, where
$$a=\sum_{j=1}^ra_jc_j\in
\bigoplus_iV_{ii}. $$
We set $\lambda_j$ the linear form on $\cal A$ given by
$\lambda_j(H)=a_j$.
It is clear that the Peirce decomposition gives also
 a simultaneous diagonalization
of all $H\in \cal A$,
namely
\begin{equation}
    \qquad\qquad Hx=L(a)x=\frac{\lambda_{i}(H)+\lambda_{j}(H)}{2}x
\hspace{3cm}  x\in V_{ij}\,.
    \label{diag}
\end{equation}
Let $A=\exp \cal A$. Then $A$ is an abelian group, and this is the
Abelian group  in the Iwasawa decomposition of $G$. We  now
describe the nilpotent part $N_{0}$. Its Lie algebra $\cal N_{0}$
is the space  of elements $X\in{\cal G}$ such that, for all $i\leq
j$, $$XV_{ij}\subset \bigoplus_{k\geq l\;;\;( k,l)>
(i,j)}V_{kl},$$
 where the pairs ordered lexicographically. Once
$\cal N_{0}$ is defined, we define $\cal S_{0}$ as the direct sum
$\cal N_{0}\oplus \cal A$. The groups $S_{0}$ and $N_{0}$ are then
obtained by taking the exponentials. It follows from the
definition of $\cal N_{0}$ that the matrices of   elements of
$\cal S_{0}$ and  $S_{0}$, in the orthonormal basis
$\{e^\alpha_{ij}\}$, are upper-triangular.

The solvable group $S_{0}$ acts simply transitively on $\Omega$.
This
may be found in \cite{FK} Chapter VI, as well as the precise
 description of $\cal N_{0}$ which will be needed later. One has
$$
{\cal N}_0=\bigoplus_{i<j\leq r}{\cal N}_{ij}, $$
where
$$
\cal N_{ij}=
\{z\qed c_i\ :z \ \in V_{ij}\}.
$$
This decomposition corresponds to a diagonalization of the adjoint
action of ${\cal A}$
since
\begin{equation}
    \qquad\qquad [H,X]=\frac{\lambda_j(H)-\lambda_i(H)}{2}X , \ \
\ X\in \cal N _{ij}\,.
    \label{comN}
\end{equation}
Finally, let $V^{\C}= V+iV$ be the complexification of $V$.  We
extend
the action of $G$ to $V^{\C}$ in the obvious way.

\subsection {Preliminaries on irreducible symmetric Siegel domains
of
type II}
We consider the Siegel domain defined by an irreducible symmetric
cone $\Omega$ and an additional complex vector space
${\cal Z}$
together with a Hermitian symmetric bilinear mapping
$$
\Phi:\ {\cal Z}\times{\cal Z}\to V^{\C},
$$
such that
\begin{align*}
&\Phi(\z,\z)\in \ov{\Omega},\ \
\z\in{\cal Z}\, ,\\
& \Phi(\z,\z)=0\ \text{ implies }
\z=0.
\end{align*}
The Siegel domain  associated with these data is defined
as
\begin{equation}
    {\cal D}=\{(\zeta,z)\in{\cal Z}\times V^{\C}:\
\Im z-\Phi(\zeta,\zeta)\in\Omega\}.
    \label{Siegel}
\end{equation}
It is called of tube type if $\cal Z$ is reduced to $\{0\}$.
Otherwise, it is called of type II. There is a representation
$\s: S_0\ni s\mapsto \s(s)\in GL({\cal Z})$
such that
\begin{equation}
s\Phi(\z, w)=\Phi(\sigma(s)\z,
\sigma(s)w)\, ,\label{repr}
 \end{equation}
and such that all automorphisms $\sigma (s)$, for $s \in A$, admit
a joint
diagonalization (see \cite {KW}). To reduce notations, we shall as
well denote by $\s $
the corresponding representation of
the algebra $\cal S_0$. For $X\in \cal S_{0}$, (\ref{repr}) implies
that
\begin{equation}
X\Phi(\z,w)=\Phi(\sigma(X)\z,w)+
\Phi(\z,\sigma(X)w)\, .\label{repr-alg}
 \end{equation}
As an easy consequence, one can prove that the only possible
eigenvalues for $\s(H)$, with $H\in \cal A$, are
${\lambda_j}(H)/{2}$, for $j=1,...,r$. So we may write $$  {\cal
Z} = \bigoplus _{j=1}^r{\cal Z} _j $$ with the property that
\begin{equation}
    \s (H)\zeta =
    \frac{\lambda_j(H)}{2}\zeta, \ \ \ \ \zeta \in
{\cal Z} _j.
    \label{zi}
\end{equation}
Moreover, all the spaces ${\cal Z} _j$ have the same
dimension\footnote{
In fact, the present study generalizes to all homogeneous Siegel
domains related to irreducible symmetric cones for which this last
property is satisfied.}. A proof of
these two facts may be found in \cite{DHMP}. We call $\chi$ the
dimension of ${\cal Z} _j$ for $ j=1,...,r$. Let us remark, using
(\ref
{repr-alg}) and (\ref{zi}), that for
 $\z ,w \in \cal Z_{j}$, we have
$L(c_{j}) \Phi(\z ,w)=\Phi(\z ,w)$. Therefore, $\Phi(\z
,w)=Q_{j}(\z ,w)c_j$, for $\z ,w \in \cal Z_{j}$. Moreover,
$\langle c_j,\Phi(\z ,\z)\rangle >0$ for $\z\in \cal Z_{j}$ and so
the Hermitian form $Q_{j}$ is positive definite on ${\cal Z} _j$.

The representation $\s$ allows to consider ${ S}_0$ as a group
of  holomorphic automorphisms of ${\cal D}$. More generally,
the elements $\zeta\in{\cal Z}$, $x\in V$ and $s\in S_0$
act on
${\cal D}$ in the following way:
\begin{align}
\z\cdot (w,z)&=(\z+w,z+2i\Phi(w,\z)
+i\Phi(\z,\z))\, ,\nonumber \\
x\cdot (w,z)&=(w,z+x)\, ,\label{mult}\\
s\cdot (w,z)&=(\sigma(s)w,sz)\, .\nonumber
\end{align}
We call $N(\Phi)$ the group corresponding to the first two actions,
that is $N(\Phi)={\cal Z}\times V$ with the product
\begin{equation}
    (\z,x)(\z',x')=(\zeta+\z',x+x'+2
\Im\,\Phi(\z,\z')).
    \label{prod}
\end{equation}
All three actions generate  a solvable Lie group
$$
S=N(\Phi)S_0=N(\Phi )N_0A=NA,
$$
which identifies with a group of holomorphic automorphisms acting
simply transitively on  ${\cal D}$.
 The group $N(\Phi)$, that is two-step nilpotent, is a
normal subgroup of $S$.
 The  Lie algebra ${\cal S}$ of $S$ admits the decomposition
 \begin{equation}
\cal S = \cal N(\F)\oplus \cal S_0
= \bigg( \bigoplus^r_{j=1}{\cal Z}_j \bigg)
\oplus \bigg( \bigoplus_{i\le j}\,V_{ij} \bigg)
\oplus \bigg( \bigoplus_{i<j}\,{\cal N}_{ij} \bigg)
\oplus\cal A.\label{direct}
\end{equation}
  Moreover, by (\ref {diag}), (\ref {comN}) and (\ref {zi}), one
knows the adjoint
  action of elements $H\in \cal A$:
  \begin{align}
[H,X]&=\frac{\lambda_j(H)}{2}X\hskip1.7cm\text{for}\quad
X\in{\cal
Z}_j,\nonumber\\
[H,X]&=\frac{\lambda_i(H)+\lambda_j(H)}{2}X\quad\text{for}\quad
X\in
V_{ij},\label{action}\\
[H,X]&=\frac{\lambda_j(H)-\lambda_i(H)}{2}X\quad\text{for}\quad
X\in{\cal
N}_{ij}.\nonumber\end{align}

  Since $S$ acts simply transitively on the domain ${\cal
 D}$, we may identify $S$ and ${\cal D}$. More precisely, we define
\begin{equation}
    \theta: S\ni s\mapsto \theta (s)=s\cdot{\bf e}\in \cal D
    \label{theta},
\end{equation}
where ${\bf e}$ is the point $(0,ie)$ in ${\cal D}$. The Lie
algebra
$\cal S$ is then identified with the tangent space of  ${\cal D}$
at
${\bf e}$ using the differential $d\th_{e}$. We  identify $e$ with
the unit element of
$S$.
We  then transport both the Bergman metric $g$ and the complex
structure $\cal J$ from $\cal D$ to $S$, where
they become left-invariant tensor fields on $S$. We still write
$\cal J$
for the complex structure on $S$. Moreover, the complexified
tangent
space $T_{\bf e}^{\C}$ is identified with the complexification of
$\cal S$,
which we denote by $\cal S^\C$. The decomposition $T_{\bf
e}^{\C}=T_{\bf
e}^{(1,0)}\oplus T_{\bf e}^{(0,1)}$ is transported into
\begin{equation}
    \cal S^\C=\sq\oplus\sp\,.
\end{equation}
 Elements of $\cal S^\C$ are identified
with left
invariant vector fields on $S$, and are called left invariant
holomorphic
vector fields when they belong to $\sq$. The conjugation operator
 exchanges $\sq$ and $\sp$, while the transported operator $\cal J$
 coincides with $i\text{Id}$ on $\sq$, and to $-i\text{Id}$ on
$\sp$. The K\"ahlerian metric given by the Bergman metric can be
seen as a Hermitian form on $\sq$, and orthonormality for left
invariant  holomorphic vector fields means orthonormality for the
corresponding elements in $\sq$.

\medskip
Now, let us construct an orthonormal
basis
of left invariant holomorphic vector fields. We first build a basis
in ${\cal
 S}$. To do this, we use the  decomposition given in
(\ref{direct})
 and give a basis for each block.

We have already fixed an orthonormal basis $\{e^\a_{jk}\}$ in $V$
corresponding to the
Peirce decomposition chosen. For $j<k$ and $1\leq \a\leq d$, we
define
$X^\a_{jk}\in V_{jk}$ and $Y^\a_{jk}\in{\cal N}_{jk}$
as the left-invariant vector fields on $S$ corresponding
to $e^{\a}_{jk}$ and $2e_{jk}^{\a}\,\qed c_j$, respectively.
For each $j$ we define $X_{j}$ and
 $H_j$ as the left-invariant vector fields on $S$ corresponding to
$c_{j}\in V_{jj}$ and $L(c_j)\in \cal A$, respectively.
It remains to choose a basis of each $ \cal Z_{j}$. We choose
for $e_{j\a}$ an orthonormal basis of $\cal Z_{j}$ related to
$4Q_{j}$, where $Q_{j}$ is the quadratic form defined above.
For $z_{j\a}=x_{j\a}+iy_{j\a}$ the corresponding coordinates, we
define
$\sx^\a_j,\sy^\a_j$ as the left-invariant vector fields on $S$
which coincide with  $\partial_{x_{j\a}}$ and $\partial_{y_{j\a}}$
at ${\bf e}$.

Finally, we define $$Z_{j}=X_{j}-iH_j ,\hspace{2cm}  Z_{jk}^\alpha
=X_{jk}^\alpha-iY_{jk}^\alpha ,\hspace{2cm}
\sz_j^\alpha=\sx_j^\alpha-i\sy_j^\alpha\,. $$

We can now state the following lemma.
\begin{lem}
    The left invariant vector fields $Z_{j}$, for $j=1,\cdots, r$,
$Z_{jk}^\alpha$, for $j< k\leq r$  and $\a=1,\cdots, d$, and
    $\sz_j^\alpha$ for  $j=1,\cdots, r$ and  $\a=1,\cdots, \chi$,
    constitute an orthonormal basis of holomorphic left invariant
    vector fields.
\end{lem}
\bpf
This lemma is already contained in \cite{DHMP}, to which we refer
for
details. To
prove that $Z_{j}$,  $Z_{jk}^\alpha$, and  $\sz_j^\alpha$ are
holomorphic vector fields, it is sufficient to prove that
$$
{\cal J}(X_{j})=H_j ,\hspace{2cm}  {\cal
J}(X_{jk}^\alpha)=Y_{jk}^\alpha ,\hspace{2cm}
{\cal J}(\sx_j^\alpha)=\sy_j^\alpha\,.
$$
To do this, we compute the image of the vector fields $X_{j}, H_j,
X_{jk}^\alpha ,Y_{ij}^\alpha ,\sx_j^\alpha$, and $\sy_j^\alpha$  by
the differential
$d\theta_e$. We find the  following tangent vectors
at {\bf e}: $\partial _{x_{jj}}$,  $\partial
_{y_{jj}}$, $\partial
_{x^{\alpha }_{jk}}$, $\partial _{y^{\alpha }_{jk}}$, $\partial
_{x_{j\alpha }}$,
and $\partial _{y_{j\alpha }}$. Here the coordinates that we have
used in
$\cal Z \times V^{\C}$
are given by
$$(\z ,z)=\Big (\sum _{j,\alpha }(x _{j\alpha}+iy _{j\alpha})
e_{j\alpha },
\sum _{i\leq j , \alpha }(x^{\alpha }_{ij}+i y^{\alpha
}_{ij})e^{\alpha }_{ij}\Big ).
$$
The assertion follows at once, using the complex structure in
$\cal Z \times V^{\C}$.

To show orthonormality,  it is possible to use  Koszul's
formula which allows to get the Bergman metric from the adjoint
action.
This is done in \cite{DHMP}, Lemma (1.18).
\epf

\subsection {Hua-Wallach systems
for irreducible symmetric Siegel domains}
We now compute the operator $\HH$ in the orthonormal basis that we
have
 built in the previous subsection. In fact, it is enough to compute
the following operators, called  {\it strongly
    diagonal} $HW$ operators, and defined by
$$\qquad\qquad \HH _j f=(\HH f\cdot Z_j, Z_j)\ ,\qquad\qquad
j=1,\cdots,r.
$$
We have the following proposition.
    \begin{prop} The strongly diagonal HW operators $\HH _j$ are
       \begin{equation}
        \HH _j= \sum _{\a}\sL _j^{\a}+2 \De _j + \sum _{k<j} \sum_{\alpha}
    \De^{\a}_{kj}+
\sum _{l>j} \sum_{\alpha}\De ^{\a}_{jl}\,,\label{Hj}
    \end{equation}
where
\begin{align}
\Delta_j&=X_{j}^2+H_j^2-H_j \,\nonumber \\
{\cal L}_j^{\a}&=(\sx_j^{\a})^2+(\sy^{\a}_j)^2-H_j\label{block}\,
\\
\Delta_{ij}^{\a}&=(X_{ij}^{\a})^2+(Y_{ij}^{\a})^2-H_j\, .\nonumber
\end{align}
\end{prop}
\bpf
We first  compute the curvature tensor $R( \ov Z, Z)$, with
$Z=Z_1,\cdots, Z_r$. From (\ref{connec}), we know that, for
$Z,W$ in $\sq$,
\begin{equation}
   \del_{\overline Z}W = \pi_{\sq}([\overline Z,W])=
\pi_{\sq}([\overline
   Z,(W+\ov W)]),
    \label{connec2}
\end{equation}
where  $\pi_{\sq}$ denotes the projection from $\cal S^\C$ onto
$\sq$. We claim that
\begin{lem}\label{nabla}
  The following identities hold:
    \begin{align*}
 \del_{\overline Z_{j}}Z_{k} =& i\d _{jk}Z_j\\
\del_{\overline Z_{j}}Z^{\a }_{kl}=&\frac{i}{2}(\d _{lj}Z^{\a
}_{kj}+\d
_{kj}Z^{\a
}_{jl})
\ \ \text{if} \quad k<l\\
\del_{\overline Z_{j}}\sz _k^{\a} =& \frac{i}{2}\d _{jk}\sz ^{\a
}_{j}.
\end{align*}
\end{lem}
\bpf
In the computation, we have seen that we may replace the three left
hand sides of the formulas above by $2\pi_{\sq}([\overline
Z_j,X_k])$,
$2\pi_{\sq}([\overline Z_j,X^{\a }_{kl}])$ and
$2\pi_{\sq}([\overline Z_j,\sx _k^{\a}])$, respectively. Moreover,
if we replace
$\overline
Z_j$ by $iH_{j}$ in these three expressions, we obtain the right
hand
sides, by virtue of (\ref{action}). Thus the lemma follows, once we
prove that all brackets
$[\overline Z_j,X_k]$, $[\overline Z_j,X^{\a }_{kl}]$, and
$[\overline Z_j,\sx _k^{\a}]$ vanish. This last fact follows from
a
standard argument. One proves that each of these vector fields is
annihilated by all endomorphisms $\text{ad} H-\lambda
(H)\text{Id}$, with
$H\in \cal A$, for a
value  $\lambda (H)$ that is not an eigenvalue of $\text{ad} H$ for
some $H$. So it
vanishes.
\epf

Let us go on with the proof of the proposition. It is easy to
deduce
the action of $\del_{Z}$ on $\sq$ from the one of $\del_{\ov Z}$.
Indeed, since the action of $S$ preserves the Hermitian scalar
product, and since $Z$ is left-invariant,
$$0=Z\cdot (U,V)=( \del_{Z}U,V)+(U, \del_{\ov Z}V)$$
for any couple $U$, $V$ of left-invariant holomorphic vector
fields.
So the endomorphism of $\sq$  defined by $\del_{Z}$ is
 the opposite of the adjoint endomorphism defined by $\del_{\ov
Z}$. It follows from the matrix representation  given in the lemma
that they are equal, and they commute. So, for $U\in\sq$, $$R(\ov
{Z_{j}}, Z_{j})U=-\del_{[\ov
Z_{j},Z_{j}]}U=-2i\del_{\ov{Z_{j}}}U$$ since $[\ov
{Z_{j}},Z_{j}]=2iX_{j}=i(Z_{j}+\ov {Z_{j}})$. Using again Lemma
\ref{nabla} and the expression of $\HH f$ given in
(\ref{Hua-inv2}), we see that $$ \HH _j= \sum _{\a}\De(\sz
_j^{\a},\sz _j^{\a})+2 \De (Z_j,Z_j) + \sum _{k<j}
\sum_{\alpha}\De(Z^{\a}_{kj}, Z^{\a}_{kj})+ \sum _{l>j}
\sum_{\alpha}\De(Z^{\a}_{jl},Z^{\a}_{jl})\,.$$ We refer to \cite
{DHMP} for the computation of $\De(\sz _j^{\a},\sz _j^{\a})$,
$\De (Z_j,Z_j)$, and $\De(Z^{\a}_{kj}, Z^{\a}_{kj})$. \epf

We also refer to \cite {DHMP} for the computation of the
Laplace-Beltrami operator $\De$,
\begin{equation}
 \De =\sum _j \De _j + \sum _{k<j}
\sum_{\a}\De^{\a}_{kj}
+\sum _{j,\a }\sL _j^{\a}.\label{L-B}
\end{equation}
It is proved in \cite{DHP} that the Laplace-Beltrami operator is
the
trace of the  operator $\HH$.

\medskip
All results, up to now, are also valid for the tube domain
$T_{\W}=V+i\W$,
 which identifies  with the subgroup $VS_0$ of the group $S$ and
 appears as a particular case. Left
 invariant differential operators act from the right. Therefore,
 we can identify
 left invariant differential operators on the tube domain with left
invariant differential operators on the domain $\sd$ itself. We
add a
 subscript or superscript for such operators coming from the tube
domain, and
 define $\HH _{j}^T$, $j=1,\cdots,r$, and $\De_{T}$ as the
operators
 coming from the strongly
 diagonal $HW$ operators for the tube domain and the
Laplace-Beltrami
 operator, respectively. Then, we
 have the following corollary, the proof of which is immediate:
 \begin{cor} The following identities hold:
   \begin{equation}
        \HH _j^T=2 \De _j + \sum _{k<j} \sum_{\alpha}\De
^{\a}_{kj}+
 \sum _{l>j} \sum_{\alpha}\De ^{\a}_{jl}\,;\label{HjT}
    \end{equation}
 \begin{equation}
\De _T=\sum _{j=1}^r\HH _j-\De\,.\label{L-B-T}
    \end{equation}
\end{cor}
\subsection{Induction procedure}
We collect in this subsection some information and some notations
which will be used in all proofs which are based on induction on
the rank of the cone. So, here, we assume that $r>1$. We first
define $$ {\cal A}^-=\text{lin} \{L(c_1),...,L(c_{r-1})\}\ \ \
\text{and}\ \
\
 {\cal A}^+=\text{lin} \{L(c_r)\}\,,
$$
and, in an analogous way,
$${\cal N}_0^-
=\bigoplus _{i<j\leq r-1} \sn _{ij}
\ \ \ \text{and}\ \
\
 {\cal N}_0^+=\bigoplus _{j=1}^{r-1} \sn _{jr}\,.
$$
${\cal N}_0^+$ is an ideal of ${\cal N}_0$, while ${\cal N}_0^-$ is a
subalgebra. Clearly ${\cal A}={\cal A}^-\oplus {\cal A}^+$ and
${\cal N}_0={\cal N}_0^-\oplus {\cal N}_0^+$.

Next, we define $A^+, A^-, N_{0}^+, N_{0}^-$ as the exponentials
of the corresponding Lie subalgebras. Then $S_{0}^-=N_{0}^-A^-$ is
the solvable group corresponding to the cone $\Omega^-$,
determined by the frame $c_{1}, \cdots, c_{r-1}$, which is of rank
$r-1$ as we wanted. The underlying space $V^-$ for $\Omega^-$ is
the subspace $$V^-=\bigoplus_{1\leq i\leq j< r}V_{ij}\,. $$ We
will make an extensive use of the fact that $$
 A=A^-A^+\ \ \text{and}\ \  N_0=N_0^-N_0^+
$$ in the sense that the mappings $$ A^-\times A^+\times
\ni(a^-,a^+)\mapsto a^-a^+\in A, $$ and $$ N_0^-\times N_0^+\ni
(y^-,y^+)\mapsto y^-y^+\in A $$ are diffeomorphisms.

Now, let us define
$$\sz ^-=\bigoplus _{j=1}^{r-1}\sz _j\,.$$
Then it is easily seen that
 $\sz ^-\times \sz ^-$ is mapped by $\Phi$ into the subspace $(V^-) ^{\C}$.
 Moreover,
$\Phi(\zeta,\zeta)$ belongs to $\Omega^-$ when $\z\in \sz ^-$.
So, we may define the Siegel domain $\sd^-$ as
$$
{\cal D}^-=\{(\zeta,z)\in{\cal Z}^-\times( V^-)^{\C}:\
\Im z-\Phi(\zeta,\zeta)\in\Omega^-\}.
$$
Let us define $\sn (\F )^-=\sz ^-\oplus V^-$ and
$$\sn (\F )^+=\sz_{r}\oplus \bigoplus_{j\leq r} V_{jr}\,.$$

Then, again, $\sn (\F )^-$ is a subalgebra and $\sn (\F )^+$ is an
ideal of
$\sn (\F )$.
We define $N(\Phi)^-$ and $N(\Phi)^+$ as their exponentials. Then
 $N(\F)$ is a semi-direct product
$$
N(\F )=N (\F )^-N (\F )^+.
$$
Clearly $N(\F )^-$ is the nilpotent step
two
``boundary" group corresponding to $\sd ^-$.

Finally, we want to decompose the group $N$. Let $N^-= N(\F
)^-(N_0)^-$, and $N^+= N(\F )^+(N_0)^+$. Then $N$ is a semi-direct
product $N=N^- N^+$. Moreover, the whole group $S$ may be written
as $$ S=N^- N^+A^-A^+=N^-A^- N^+A^+. $$ Clearly,
$
S^-=N^-A^-
$
is the solvable group acting simply transitively on $\sd^-$.

\section{ Poisson integrals}
The aim of this section is to prove the following partial result in view of
the main theorem.
\begin{thm} Let $F$ be a bounded function on $S$
annihilated by
$\De $ and by $\HH_j$, for $j=1,...,r$. Then
$$
\HH_j^TF=0\ \ \ \text{for}\ j=1,...,r,
$$
and
$$
\sL _jF=\sum _{\a}\sL _j^{\a }F=0\ \ \ \text{for}\ j=1,...,r.
$$
\label{sectP}\end{thm}
>From the formulas of the previous section, it is clear that a bounded
function on the domain $\sd$ which is $\HH$--harmonic satisfies the
assumptions. Moreover, the first statement implies the second one. To
prove the first one, we shall use the characterization of $\HH$--harmonic
functions in terms of Poisson-Szeg\"o integrals on tube domains. More
precisely, following Hua \cite{Hua} and \cite{FK}, it is sufficient to prove that
$F$, considered as a function on the tube domain $T_{\W}=V+i\W$, is
the Poisson-Szeg\"o integral of some bounded function on $V$. To do
this, our main tool will be the possibility to write $F$, in different
ways, as a Poisson integral related to some elliptic operators which
annihilate $F$.

Let us first give some notations. From the last
section, we know that every $g\in S$ may be written in a unique way as a product
$(\z, x)na$, with $(\z, x)\in N(\Phi)$ and $n\in N_{0}$. We write
$\pi$ for the projection on $N(\Phi)$, given by $\pi(g)=(\z, x)$, and
$\tilde \pi$ for the projection on $N$, given by $\tilde\pi(g)=(\z, x)n$.

We first recall previous results of two of the authors. Even if they are
valid in the more general context of a semi-direct product, we give
them in the present context. We consider elliptic operators which may
be written as
\begin{equation}
L=\sum_{j=1}^{r} \a _j\sL _j+\sum_{j=1}^{r} \b _j\HH_j^T
\label{L}\end{equation}
with $\a_{j}$ and $\b_{j}$ positive constants.
Then $L$ is a sum of square of vector fields plus a first order term
$Z=Z(L)$, which is called the drift, and may be written as
$Z=-\sum \gamma_{j}H_{j}$, with $\gamma_{j}=\a _j\x +(2+(j-1)d)\b _j+d\sum _{k<j}\b _k$.
It follows from \cite{DH} and \cite{R} that the maximal boundary of $L$
can be easily computed (it depends on the
signs of $\lambda_{j}(Z)- \lambda_{i}(Z)$ for $i<j$). In particular, it is equal to
$N(\Phi)$ if the sequence $\gamma_{j}$ is a non--increasing sequence,
and  to $N$ if it is an increasing sequence. Let us summarize
the results that we shall use in the next proposition.
\begin{prop}Let $L$ be given by {\rm(\ref{L})}, and $ \gamma_{j}$ as
above.\\
$(i)$ If $L$ is such that $\gamma_{j}$ is a non--increasing sequence,
there is a unique positive, bounded,
 smooth function
$P_{L} $ on $N(\Phi)$ with $\int _{N(\Phi)}P_{L}(y)dy=1$ such that bounded
$L$-harmonic
functions on $S$ are in one-one correspondence with
$L^{\infty}(N(\Phi))$ via
the Poisson integral
\begin{equation}
F(s)= P_{L}f(s)=\int _{N(\phi)}f(\pi (sw))P_{L}(w) \ dw.
\label{poisson}
\end{equation}
$(ii)$ If $L$ is such that $\gamma_{j}$ is an increasing
sequence, there is a unique positive, bounded,
 smooth function
$\tilde P_{L} $ on $N$ with $\int _{N}\tilde P_{L}(y)dy=1$ such that bounded
$L$-harmonic
functions on $S$ are in one-one correspondence with
$L^{\infty}(N)$ via
the Poisson integral
\begin{equation}
F(s)=\tilde P_{L}f(s)=\int _{N}f(\tilde\pi (sy))\tilde P_{L}(y) \ dy.
\label{poisson2}\end{equation}
Moreover, for each given $\et>0$, we may choose the coefficients
$\a_{j}$ and $\b_{j}$ so that $(i)$ holds, and that
\begin{equation}\int _{N(\F )}\t (w)^{\et }P_L(w)\ dw <\infty
\label{int}\end{equation}
where $\t (w)$ is the distance of $w$ from the unit element $e\in N(\F )$ with respect to any
left--invariant Riemannian metric.
\end{prop}
As we said, $(i)$ and $(ii)$ may be found in \cite{DH} and
\cite{R}. The integrability condition may be found in \cite{D},
Theorem (3.10): a sufficient condition for (\ref{int}) is that
$$\et \sum 2\b _j\lambda(H_j)^2+\lambda(Z)<0, $$ for all linear
forms on $\cal A$ of the form
$\lambda=\frac{\lambda_k+\lambda_p}{2}, \frac{\lambda_k}{2}$. The
fact that this condition may be satisfied is elementary.

We have chosen to add a tilde every time that we are concerned with
an operator whose maximal boundary is the whole group $N$. We then
define  $P_{L}$ as an integral,
\begin{equation}
P_{L}(w)=\int _{N_0}\tilde P_{L}(wy)\ dy.
\label{intL}\end{equation}
Let us remark that, in this case, the functions $F$ which may be written as
\begin{equation}
F(s)= P_{L}f(s)=\int _{N(\phi)}f(\pi (sw))P_{L}(w) \ dw,
\label{poisson3}\end{equation}
with $f$ a bounded function on $N(\Phi)$, constitute a proper
subspace of the space of bounded functions which are annihilated by
$L$. It is in particular the case for the Laplace-Beltrami operator,
which is obtained for the values $\a=2\b=1$, and has maximal
boundary $N$.

\medskip

The main step in this section is the next proposition.
It has been proved in \cite{DHP} for general homogeneous Siegel domains
(non
necessarily symmetric), and for more general operators.
However, in the case of symmetric Siegel domains, which is the case
under consideration, the
proof
 may be simplified considerably.
 We include it for
the reader's
convenience.
\begin{prop}
   Let $F$ be a bounded function on $S$
annihilated by
$\De $ and by $\HH_j$, for $j=1,...,r$. Then, there exists a bounded
function $f$ on $N(\Phi)$ such that $F$ may be written as
$$F(s)= P_{\De}f(s)=\int _{N(\phi)}f(\pi (sw))P_{\De}(w) \ dw.$$
\label{disting}\end{prop}
\bpf
We already know that there exists some bounded function $\tilde f$ on $N$
such that $F$ may be written as $\tilde P_{\De}\tilde f$. Moreover, we may
assume that $\tilde f$ is a continuous function and prove that, in this case,
$f$ is the restriction of $\tilde f$ on $N(\phi)$. Indeed, in the general
case, we consider
the sequence of functions $F_{m}$ defined by
$$
F_m(s)=\int _N\f _m(n)F(n^{-1}s)\ dn=\tilde P_{\De}(\f _m*\tilde f)(s).
$$
 with $\f _m$ an approximate identity
which is compactly supported and of class ${\cal C}^\infty$.
Clearly $F_{m}$ tends to $F$ pointwise. Let us assume that we have
already proved the proposition for continuous functions. Then
$F_m=P_{\De}(f_m)$. All the functions $(f_m)$ are bounded by $\|
f\| _{L^{\infty}}$, so we can extract a $^*$--weak convergent
sequence which converges to $f$. Then $P_{\De}(f_m)$ converges to
$P_{\De}(f)$ pointwise. Hence $F=P_{\De}f$.

So, let $\tilde f$ be a bounded continuous function on $N$, and
let $F=\tilde P_{\De}\tilde f$. To prove the proposition, we want
to prove that, for each fixed $w\in N(\Phi)$, the function
$y\mapsto \tilde f(wy)$ is constant on $N_{0}$. Indeed, assume
that it is the case and  denote by $f$ the restriction of $\tilde
f$ to $N(\Phi)$. Then, for  $s=wya$ with $w\in N(\Phi)$, $n\in
N_{0}$ and $a\in A$, we can write
\begin{align*}
F(wya)=\tilde P_\De \tilde f(wya)=&\int _{N(\F )N_0}\tilde
f(wyavua^{-1})\tilde P_\De(vu)\ dvdu,\\ =&\int _{N(\F
)N_0}f(wyava^{-1}y^{-1})\tilde P_\De(vu)\ dvdu,\\ =&P_\De f(wya).
\end{align*}
Let us finally remark that  it is sufficient to prove that $y\mapsto \tilde f(y)$ is
constant on $N_{0}$. Indeed, once we have proved this, for each $w\in N(\Phi)$
we have the same conclusion with $F$ replaced by $_wF$, with  $_wF(g)=F(wg)
=\tilde P_\De(_w\tilde f)(g)$. Again $_w\tilde f(y)=\tilde f(wy)$ is
constant, which we wanted to prove.

So, let us show that  $y\mapsto \tilde f(y)$ is
constant on $N_{0}$. Let us define
\begin{equation}
F_H(wya)=\int _{N_0}\tilde f(yaua^{-1})\left (\int _{N(\F )}\tilde P_\De(vu)\
dv\right)\ du.\label{FH}\end{equation}
We claim that
 \begin{equation}
     F_H(g)=\lim _{t\to -\infty}F((\exp tH)g), \label{lim}
 \end{equation}
 where $H$ is the vector field
 $H=\sum _{j=1}^rH_j$. Indeed,  writing
$$
F(g)=\int _{N(\F )N_0}\tilde f(\tilde\pi (gvu))\tilde P_\De(vu)\ dvdu,
$$
we have
$$
F((\exp tH)wya)=F(w_ty_ta\exp tH)=\int _{N(\F )N_0} \tilde f(w_ty_tav_tu_ta^{-1})
\tilde P_\De (vu)\ dvdu.$$
For an element $g$ of $N$ we
have used the notation
$$
g_t=(\exp tH)g(\exp (-tH)).
$$
It follows from (\ref{action}) that $u_t=u$ for every
$u\in N_0$, and that $w_{t}$ tends to the unit element. This implies
(\ref{lim}).
We now claim that
\begin{align}
F_H(wya)=&F_H(ya) \label{FH1}\\ \HH _jF_H=&0\ \ \text{for}\
j=1,...,r\label{FH2} \ \ \text{and} \ \De F_H=0\\ HF_H=&0.
\label{FH3}
\end{align}
We have already proved  (\ref{FH1}). Then (\ref{FH2}) follows from
the fact that left and right translations commute. So, for every
$t$, $F((\exp tH)g)$ is annihilated by the $HW$ operators and the
Laplacian. To see (\ref{FH3}), we use again the fact  that $u_t=u$
for every $u\in N_0$ and the formula (\ref{FH}) to obtain that $$
F_H(ya(\exp tH))= F_H(ya).$$ Then (\ref{FH3}) follows at once.

Finally, uniqueness in Proposition \ref{int} implies that $y\mapsto
\tilde f (y)$ if and only if $F_{H}$ is constant. To prove that $F_{H}$ is constant,
we consider the function $G$ defined on $N_{0}A$ by $G(ya)=
F_{H}(ya)$. Then clearly $G$ is annihilated by all operators
\begin{equation}
    D_j=-\x H_j+2(H^2_j-H_j)+\sum _{i<j}\sum_{\a}((Y^\a_{ij})^2-H_j)+
\sum _{j<k\leq r}\sum_{\a}((Y^\a_{jk} )^2-H_k)
\label{Dj}\end{equation}
and by $H$. So, to complete the proof, it is sufficient to prove the
following lemma.
\begin{lem} Let $G$ be a bounded function on $N_0A$ which annihilated
    by the operators $H, D_1,..., D_r$. Then $G$ is constant.
    \end{lem}
    \bpf
There is nothing to prove when $r=1$. For $r=2$, let us remark that
$G$, which is annihilated by $H_{1}+H_{2}$, is also
annihilated by $H_{1}-H_{2}$ since
$$
(D_1-D_2)G=-(\x +2)(H_1-H_2)G=0.
$$
Therefore, $H_1G=0$ and $H_2G=0$ and so, $G$ is a bounded function
on the
Abelian group $N_0=\R ^d$ annihilated by the Laplace operator.
Hence $G$ is
constant.

Let us now consider $r>2$. We assume that the lemma has been proved
with $r$ replaced by $r-1$. We write $G$ as a Poisson integral with respect to the
operator
$$
D=\sum \a _j D_j, $$
for which the drift $Z(D)$ is equal to $Z=-\sum \gamma_{j}H_{j}$, with
 $\gamma_{j}=(\x +2+(j-1)d)\a _j+d\sum _{k<j}\a _k$.  We first remark that we may choose the
coefficients $\a_{j}$ so that the $\gamma_{j}$  decrease for
$j\geq 2$ (when $\x +2>d$, one can even find a sequence $\a_{j}$
such that $\gamma_{j}$ is decreasing, and conclude directly since
every $D$--harmonic bounded function is constant). With this
choice,  the maximal boundary of $D$ is the group $N_{1}=\exp
{\cal N}_{1}$, with $$\sn ^1=\oplus _{j>1}\sn _{1j}.$$ We also
define  $N^1=\exp {\cal N}^{1}$, with $\sn ^{1}=\oplus _{1<i<j\leq
r}\sn _{ij}$. Every $y$ in $N_{0}$ can be written in a unique way
as $y_{1}y'$, with $y_{1}\in N_{1}$ and $y'\in N^{1}$. We define
$\pi_{1}$ by $\pi_{1}(ya)= y_{1}$. Then, (see \cite {DH}), there
exists   functions $\nu_{D}$ and $\phi$ such that $$ G(ya)=\int
_{N_1}\phi (\pi_{1}(yau))\n _D(u)\ du. $$ The function $\phi$ is
bounded, and we can assume as before that it is continuous. Using
notations of the subsection 2.5. on the induction procedure, we
can also write $y\in N_{0}$ as $y^+y^-$. When $y$ is in $N_{1}$,
then $y^+$ belongs to $N_{1r}=\exp \sn _{1r}$. We shall prove that
$\phi (y)$ depends only on $y^+$.  Again, to prove this, it is
sufficient to prove that $\phi (y^-)$ is  constant. Indeed, once
we have proved this, we may apply it to $_{y^+}\phi$ (with
$_{y^+}\phi(n)=\phi (y^+n))$, using the function $_{y^+}G$ in
place of $G$.

In order to prove that  $\phi (y^-)$ is  constant, let us define, as before,
\begin{align*}
   G^\#(ya)=  & \lim_{t\to -\infty}G((\exp tH_{r})ya) \\
   =&\lim _{t\to -\infty}\int _{N_1}\phi(\pi_{1}(y^{-}(y^+)_t au^{-}(u^+)_t a^{-1}))
\n _D(u^+u^{-})\ du^+du^{-} \\
 =&\int _{N_1}\phi(\pi_{1}(y^{-} au^{-} a^{-1}))\n _D(u^+u^{-})\
du^+du^{-}.
\end{align*}
Here $u_{t}=(\exp tH_r)u(\exp(-tH_r))$. We have used the fact
that $(y^{-})_{t}=y^{-}$, and $(y^+)_t$ tends to the unit element.
We have
$$
 G^\#(y^{-}y^+a^{-}a^+) = G^\#(y^{-}a^{-})= G^\#(y^+y^{-}a^{-}a^+)$$
\begin{align*}
    D^\#_j G^\#=&0, \ \ \text{for}\ j=1,...,r-1,\\
H^\# G^\#=&0,
\end{align*}
where $H^\#= H_{1}+H_{2}+ \cdots +H_{r-1}$, and, for $j=1, ...,r-1$,
$$
D^\#_j=2H^2_j-(\x+2)H_j+\sum _{i<j}((Y^\a_{ij})^2-H_j)+
\sum _{j<k\leq r-1}\sum_{\a}((Y^\a_{jk})^2-H_k).
$$
 From the induction hypothesis,  we conclude that
 $G^\#$ is constant. So $\phi(y^-)$ is also constant.  Hence
$\phi(y)=\phi(y^+)$, and, using obvious notations, we conclude that
$G$ may in fact be written as
$$G(ya)=\int _{N_{1r}}\phi (\pi_{1r}(yau))\n _D(u)\ du.
$$
Since $\exp tH_{j}$ commutes with elements of $ N_{1r}$ for $j=2,
\cdots r-1$, we conclude that $H_{j}G=0$. So $(H_{1}+H_{r})G=0$. Moreover,
\begin{align*}
   D_1G&=\Big (2H^2_1-(\x+2)H_1+\sum _{ \a }((Y^{\a }_{1r})^2-H_r)\Big )G=0\\
   D_rG&=\Big (2H^2_r-(\x+2+(r-2)d)H_r+\sum _{ \a }((Y^{\a
   }_{1r})^2-H_r)\Big )G=0,
\end{align*}
which, as in the case $r=2$, implies that $G$ is constant.
\epf

Once we have concluded for the lemma, we  conclude for the
proposition \ref{disting}.
\epf

Our next step is the following theorem.
\begin {thm} Let $f$ be a bounded  function on
$N(\F )$
and let $F=P_{\De}f$. Assume that
$$
\De _TF=0.
$$
Then
\begin{equation}\label{PS}
   F((\z ,x)ya)=\int _V f_{\z }(xyava^{-1}y^{-1})p(v)\ dv,
\end{equation}
where  $f_{\z }(x)=f(\z ,x)$ and $p$ is the Poisson - Szeg\"o
kernel for the tube domain $V+i\W$.\label{tube}
\end{thm}
\bpf
Using the same kind of proof as in the last proposition,
 we may assume that $f$ is continuous.
The maximal boundary for $\De _T$ considered as an operator
on $VS_0$
is $VN_0$. Let $\tilde p $ be the corresponding kernel on $VN_0$.
Then the function $F_{\z }$, which is defined  for
$\z \in \sz$ fixed by  $F((\z ,x)s)=F_{\z }(xs)$, may be written as
\begin{equation}
    F_\z(xya)=\int _{VN_0}g_\z(xyavua^{-1}))\tilde p (vu)\ dvdu,
\label{P-Sz}\end{equation}
where $v\in V$, $u\in N_0$.

We have also $$
F_\z(xya)=P_{\De }f((\z ,x)ya)=
\int _{N(\F )}f((\z ,x)ya(\et ,v)a^{-1}y^{-1}) P_{\De }(\et ,v)\
d\et dv.
$$
Let $a_t=\exp t(\sum _{j=1}^rjH_j)$. Then, on one hand,
$$
\lim _{t\to -\infty}F_{\z}(xya_t)=g_{\z }(xy)
$$
in *-weak topology on $L^\infty (VN_0)$ and on the other,
$$
\lim _{t\to -\infty}F_{\z}(xya_t)=f(\z ,x),
$$
pointwise.
Hence $g_{\z }(xy)=f(\z ,x)=f_{\z}(x)$. Therefore,
\begin{align*}
F_{\z}(xya)=&\int _{VN_0}f_{\z}(xyav a^{-1}aua^{-1})\tilde p (vu)\
dvdu,\\
=&\int _{V}f_{\z}(xyav a^{-1}y^{-1})(\int _{N_0} \tilde p (vu)\ du
)\ dv\,.
\end{align*}
It remains to prove that it is also equal to the right hand side
of (\ref{PS}). But this last expression is a $\Delta_T$--harmonic
function (since the Poisson--Szgeg\"o kernel is
$\Delta_T$--harmonic), with the same boundary values on $VN_0$.
This proves (\ref{PS}). \epf

\bigskip

{\sl Proof of Theorem \ref{sectP}}. Again, using Proposition
\ref{disting}, we may assume that $F=P_{\De}f$. Moreover, we may
assume that $f$ is continuous, using the same trick as in the
proof of Proposition \ref{disting}. So, it follows from Theorem
\ref{tube} that $F_{\z}$ is a Poisson-Szeg\"o integral on the tube
domain. We know  from \cite{Hua} and \cite {JK} that
Poisson-Szeg\"o integrals on symmetric tube domains are
annihilated by Hua operators, i.e. in our situation $F_{\z}$ is
annihilated by $\HH _j^T$. This finishes the proof. \epf

 \section{The proof of pluriharmonicity}

In this section we prove the following statement, which implies the main  theorem:
\begin{thm} Assume that
\begin{equation}
    \sup _{s\in S_0}\int _{N(\Phi)}|F((\z ,x)s|^2 \ d\z dx <\infty
\label{H2*}\end{equation}
and
\begin{equation}
    \De F=\HH _1F=...=\HH _rF=0.
\label{equ}\end{equation}
Then $F$ is pluriharmonic.
\label{fin}\end{thm}

We first claim that the results of the last section on bounded
functions apply to $(H^2)$ growth conditions. Indeed, we have the
following lemma. Here $L$ is an elliptic operator as in (\ref
{L}), $$ L=\sum_{j=1}^{r} \a _j\sL _j+\sum_{j=1}^{r} \b
_j\HH_j^T\,, $$ with coefficients chosen so that it has maximal
boundary $N(\Phi)$ and satisfies the integrability condition
(\ref{int}) for some $\eta$ to be chosen later.
\begin{lem}
    A function $F$ which satisfies {\rm(\ref{H2*})} and {\rm(\ref{equ})} may be
written as a
 Poisson integral
\begin{equation}
F(g)=\int _{N(\F )}f(\pi (gw))P_L(w)\ dw, \ \ f\in L^2(N(\F )), \
g\in S. \label{poisson2*}\end{equation}
\end{lem}
\bpf
We reduce to bounded functions by convolving $F$ from
the left. More precisely, let $\f _n\in C_c^{\infty}(N(\F ))$
be an approximate identity, and let
$$
F_n(g)=\int _{\nfi }\f _n(w)F(w^{-1}g)\ dw.
$$
Then $F_n$ is bounded,  and satisfies  (\ref{equ}). So it follows from
the last section that
$$
F_n(g)=\int _{\nfi }f_n(\pi (gw))P_L(w) \ dw,
$$
for an $f_n\in L^\infty (\nfi )$. Moreover,  $f_{n}$,
which may be obtained as a $^*$--weak limit in $L^\infty$, when
$t\rightarrow -\infty$,
of $F_{n}(\cdot \exp {tH})$ as well as a weak limit in $L^2(\nfi)$, is
uniformly in $L^2$.
Hence,
$$
\| f_n \| _{L^2(\nfi )}\leq \sup _s \int _{\nfi }|F(ws)|^2\ dw
$$
We may take for $f\in L^2(\nfi )$ the weak limit of a subsequence,
and get (\ref{poisson2*}).
This concludes the proof of the lemma.
\epf

To prove Theorem \ref {fin}, we may assume that $F=P_{L}f$ as
above. Moreover, eventually convolving $f$ in the group $N(\Phi)$
with a ${\cal C}^\infty$ compactly supported function as in the
last section, we may assume that\\
{\sl Assumption on $f$: it may
be written as $\phi*\tilde f$ where $\phi$ is a ${\cal C}^\infty$
compactly supported function.}

\medskip
At this point, our main tool will be harmonic analysis of the
nilpotent group $N(\F)$. Once we have
proved that  the Fourier transform of $f$ vanishes outside $\Omega\cup
-\Omega$, one concludes easily like in \cite {DHMP}.

\bigskip

Let us first recall  some basic facts about Fourier analysis on
$N(\F)$, following \cite {OV}. Let $(\ ,\ )$ be the Hermitian
scalar product on  $\cal Z$ for which the basis $e_{j\a}$, which
was introduced in subsection 2.3,  is orthonormal. It coincides
with $4Q_{j}$ on each  $\cal Z_{j}$,  and these subspaces are
pairwise orthogonal. For each $\lambda \in V$, let us define the
Hermitian  transformation $M_{\lambda}:\sz \to \sz$ by $$ 4\langle
\lambda,\Phi(\z ,\omega)\rangle=(M_{\lambda}\z,\omega), \ \z
,\omega \in \sz\,. $$
 and  consider the set
$$ {\L}=\{\lambda\in V:\det M _{\lambda}\neq 0\} $$ for which the
above Hermitian form is non degenerate. Remark that it is in
particular the case for  $\lambda\in \Omega $ since we assumed
that $\Phi (\z, \z)$ belongs to $\overline \Omega\setminus \{0\}$
for all $\z\neq 0$.  The same is valid for $\lambda\in -\Omega$.
So $\det M_{\lambda}$, which is a polynomial of $\lambda$, does
not vanish identically, and ${\L}$ is an open set of full measure.
It carries the Plancherel measure (see  \cite {OV}), given by $$
\r (\lambda) d\lambda=|\det M_{\lambda}|d\lambda. $$
 Let us describe the Fock representation associated
to $\lambda\in \L$. For every $\lambda \in \L$ we define a complex
structure $\sj_{\lambda}$, which determines the representation
space $\sh _{\lambda}$. Let $|M_{\lambda}|$ be the positive
Hermitian transformation such that
$|M_{\lambda}|^2=M_{\lambda}^2$. Then $$ \sj
_{\lambda}=i|M_{\lambda}|^{-1}M_{\lambda}. $$ If $\lambda \in \W$
then $\sj _{\lambda}=iI=\sj$ coincides with the ordinary complex
structure in $\sz$. For general $\lambda$, the complex structure
$\sj _{\lambda}$ has a nice description in an appropriate basis.
Namely, there is a $\lambda$-measurable choice of an $(\ ,\ )$
orthogonal basis $e^{\lambda}_1,...,e^{\lambda}_m$ such that $$
H_{\lambda}(e^{\lambda}_j,e^{\lambda}_k)=\s _j\d _{jk} $$ with $\s
_j= \pm 1$ (depending on $\lambda$ and locally constant). In the
basis $e^{\lambda}_1,...,e^{\lambda}_m, \sj e^{\lambda}_1,...$,
$\sj e^{\lambda}_m$ of $\sz$ over $\R$ we have $$ \sj
_{\lambda}(e^{\lambda}_j)=\s _j(\sj e^{\lambda}_j) \ \text{and}\
\sj _{\lambda}(\sj e^{\lambda}_j)=-\s _j e^{\lambda}_j\,. $$ Let
$$ B_{\lambda}=\Im H_{\lambda}\,. $$
 A direct calculation shows that
 $$ B_{\lambda}(\sj
_{\lambda}e^{\lambda}_j,e^{\lambda}_k)=\d _{jk} $$ and so $$
B_{\lambda}(\sj _{\lambda}\z,\z)>0\ \text{if}\ \z \neq 0. $$ We
define $\sh _{\lambda}$ as the set of all $C^{\infty}$ functions
$F$ on $\sz$ which are holomorphic with respect to the complex
structure $\sj _{\lambda}$ and such that $$ F(\cdot ) \r
(\lambda)^{\frac{1}{2}} e^{-\frac{\pi}{2}B_{\lambda}(\sj
_{\lambda}\cdot , \cdot )}\in L^2(\sz ,dz). $$ Here $dz$ is the
Lebesgue measure related to the scalar product $(\cdot, \cdot)$ on
$\sz$.

The space $\sh _{\lambda}$ is a Hilbert space for the scalar
product $$ (F_1, F_2)_{\lambda}=\int _{\sz}F_1(\z )\bar{F}_2(\z
)e^{-\pi B_{\lambda}(\sj _{\lambda}\z , \z )}\ \r(\lambda)d\z . $$
The Fock representation $U^{\lambda}$, which is a unitary and
irreducible representation on $\sh _{\lambda}$,  is given by
\begin{equation}
U^{\lambda}(\z ,x)F(\omega )=e^{-2\pi i\<\lambda
,x\>-\frac{\pi}{2}|\z |^2+\pi \omega \bar{\z }}  F(\omega -\z )\,,
\label{repr2}\end{equation} with $\omega \bar{\z}=B_{\lambda}(\sj
_{\lambda}\omega ,\z )+iB_{\lambda}(\omega ,\z )$ and $|\z |^2= \z
\bar{\z }$. Then the Fourier transform of $f\in L^1(\nfi )$, which
we note $U^{\lambda}_{f}$, is defined as the operator on $\sh
_{\lambda}$ given by $$(U^{\lambda}_{f}F,G)_{\lambda}=\int _{\nfi
}f(\z ,x )(U^{\lambda}_{(\z ,x )}F,G)_{\lambda} dx.$$
 If $f\in L^1(\nfi )\cap L^2(\nfi )$, then the
Plancherel theorem says that $$ \int _V\| U^{\lambda }_f\| ^2
_{HS} \ \rho (\lambda ) d \lambda = \| f \| ^2 _{L^2(\nfi )}. $$
 It follows that,
for  $f\in L^2(\nfi )$, $U^{\lambda }_f$ is defined for almost
every $\lambda $ and  is a Hilbert-Schmidt operator.

Now we write an orthonormal basis of $\sh _{\lambda} $, which
changes measurably with $\lambda$. For $\zeta\in\sz$, we note
$\zeta_{j, \lambda}$ its coordinates in the basis $e_{j}^\lambda$,
so that, in particular, $$B_{\lambda}(\sj _{\lambda} \z,
\z)=\sum_{j}|\zeta_{j, \lambda}|^2\;.$$

Given a multi-index $\a = (\a _1,...,\a _m)$, let $$ \c
^{\lambda}_{\a}=\frac{\pi ^{\frac{|\a |}{2}}}{\sqrt{\a !}}
\prod_j\z_{j,\lambda}^{\a_{j}\frac
{(1+\sigma_{j})}{2}}\bar\z_{j,\lambda}^{\a_{j} \frac
{(1-\sigma_{j})}{2}}, $$ Then every $\c ^{\lambda}_{\a} $ is
holomorphic with respect to the complex structure $\sj _{\lambda}$
and the family $\{\c ^{\lambda}_{\a}\} $ forms a $(\ ,\
)_{\lambda}$ - orthonormal basis. Indeed, one may verify that $$
(\c ^{\lambda}_{\a}, \c ^{\lambda}_{\b })_{\lambda }=
 \frac{\pi ^{\frac{|\a |+|\b|}{2}}}{\sqrt{\a !\b !}}
\prod _j\int _{\C }u ^{\a _j\frac {(1+\sigma_{j})}{2}}
\bar u ^{\a _j\frac {(1-\sigma_{j})}{2}}\bar u ^{\b _j\frac {(1+\sigma_{j})}{2}}
u ^{\b _j\frac {(1-\sigma_{j})}{2}}
 e^{-\pi  |u|^2}\ du.
$$ We finally define, for $f\in L^2(N(\Phi)$ and almost every
$\lambda$,
\begin{equation}
\hat f(\lambda ,\a ,\b )=(U^{\lambda}_f\xi _{\a }^{\lambda }, \xi
_{\b }^{\lambda })\,.\label{fourier}
\end{equation}

We may now give the main step of the proof.
\begin{lem}
Let $F=P_{L}f$ a function which satisfies the assumptions of
Theorem {\rm\ref{fin}}, with $L$ and $f\in L^2(N(\Phi))$ chosen as
above. Then, for almost every $\lambda$ and for all $\alpha$,
$\beta$, we have
\begin{equation}
    \hat f(\lambda ,\a ,\b )=0 \ \ \ \text{for}  \ \  \lambda \notin \bar \W \cup -\bar \W .
\label{four}\end{equation}
\label{four1}\end{lem}
{\sl Proof of Theorem  \ref{fin}.} For the moment, we take the lemma for granted and finish the proof of
Theorem \ref{fin}. Let us first give some notations. For $s\in
S_{0}$, we note $F_{s}$ the function defined on $N(\Phi)$ by
\begin{equation}
    F_{s}(\z, x)=F((\z, x)s)
    \label{Fs}\end{equation}
    and $ \hat F(\lambda ,\a ,\b, s )$ its Fourier transform.
   We claim that
\begin{align}
\hat F(\lambda ,\a ,\b ,s)=& e^{-2\pi \langle \lambda ,s\cdot e\>
} (U^{\lambda }_f\xi _{\a }, \xi _{\b }),\quad \text{for a.e} \
\lambda \in \bar \W ,\notag\\ =&e^{2\pi \< \lambda ,s\cdot e\> }
(U^{\lambda }_f\xi _{\a }, \xi _{\b }),\quad \text{for a.e} \
\lambda \in -\bar \W,\label{four2}\\ =& 0,\quad \quad \quad \quad
\quad \quad \quad \quad \text{for a.e}\ \lambda
 \notin \bar \W
\cup - \bar\W . \notag\end{align} Indeed, we know from Theorem
\ref{tube} that $F$ may be written as a Poisson-Szeg\"o integral,
i.e. $$ F((\z ,x)s)=\int _V f_\z(xsvs^{-1})p(v)\ dv=\int _V f(\z
,x-u)p_s(u)\ du\,, $$ with $p_s$ defined by
 $$p_s(u)=\det( s^{-1})\,p(s^{-1}\cdot u)\,.$$
 Here the element $s^{-1}$ is considered as acting on $V$.
  If $f\in L^1(\nfi)\cap L^2(\nfi )$, then

\begin{align*}
(U^{\lambda }_{F_s}\xi _{\a },\xi _{\b })=&\int _{\nfi }\int
_Vf(\z ,x-u)p_s(u) (U^{\lambda }_{(\z ,x)}\xi _{\a },\xi _{\b })\
dud\z dx\\ =& \int _V(U^{\lambda }_{f}U^{\lambda }_{(0,u)}\xi _{\a
},\xi _{\b }) p_s(u)\ du\\ =&(U^{\lambda }_{f}\xi _{\a },\xi _{\b
})\int _V e^{-2\pi i\< \lambda ,u\> }p_s(u)\ du.
\end{align*}
These formulas are still valid for a general function  $f\in
L^2(\nfi )$: only use an approximation of $f$ and the Plancherel theorem.

It remains to calculate the Fourier transform of $p_s$ for
$\lambda \in \bar \W \cup -\bar\W $. We shall do this for $\lambda
\in \bar \W$. For $\lambda \in -\bar \W$ the proof is analogous.
If $\lambda \in \bar \W$ we consider the bounded holomorphic
function on $V+i\W$ given by $$ G(z)=e^{2\pi i \< \lambda
,z\>}=e^{2\pi i \< \lambda ,x+is\cdot e\>}= e^{2\pi i \< \lambda
,x\>-2\pi \< \lambda ,s\cdot e\>}. $$ Then $G$ is the Poisson
integral of its boundary value, i.e. $$ G(z)=\int _Ve^{2\pi i \<
\lambda ,x-u\>}p_s(u)\ du. $$ Therefore,
 $$ G(is\cdot
e)=e^{-2\pi  \< \lambda ,s\cdot e\>}=\int _Ve^{-2\pi i \< \lambda
,u\>}p_s(u)\ du\,. $$ Finally, for $\lambda \in \bar \W $, we have
$$ (U^{\lambda }_{F_s}\xi _{\a },\xi _{\b })=e^{-2\pi  \< \lambda
,s\cdot e\>} (U^{\lambda }_{f}\xi _{\a },\xi _{\b }). $$ From
(\ref{four2}), a direct computation (see \cite{DHMP} for the
details) shows that $\De _jF=0$ for $j=1,...,r$. Moreover, we
already know that  $\sL _jF=0$. Then it follows from Theorem 3.1
in \cite{DHMP} that $F$ is the real part of an $H^2$ holomorphic
function. \epf

{\sl Proof of Lemma \ref{four1}.} It remains to prove the lemma. Let us remark that there is nothing to
prove for $r=1$. So the theorem is completely proved in this case. For
$r>1$, we can make the assumption that the theorem is valid for $r-1$,
and prove the lemma with this additional induction hypothesis.

We  use again the notations of the subsection 2.5. for the
induction procedure. An element $a\in A$ will be written as $a=a'a^+$, $a'\in A^{-}$,
$a'\in A^+$.
 We call $S_{0}'$ the group $N_{0}A^-$, and  $S'$ the group $NA^-$. For
 $s\in S_{0}$, we may write
 $s=ya=ya'a^+=s'a^+$.

We define a new function $F'$ on $S'$ by a limit process. More
precisely, for $(\z ,x)s'\in S'$, we define
\begin{equation}
   F'((\z ,x)s')=F'_{s'}(\z ,x)=\lim _{t\to -\infty}F((\z ,x)s'\exp t
H_{r}). \label{P'}\end{equation} Using the same arguments as
before, as well as  our assumptions on the boundary value $f$ of
$F$, one can see that this limit exists and is given by $$
F'_{s'}(\z ,x) = \int _{\nfi^{-}}f((\z ,x)s'w^{-}(s')^{-1})
P'_L(w^-) \ dw^-, $$ where $$ P'_L(w^-)=\int _{\nfi
^+}P_L(w^-w^+)\ dw^+. $$ We are now able to give a sketch of the
proof. The function $f$ may be seen as the boundary value of $F'$.
So, we will consider the Fourier transform of $F'_{s'}$. Using the
induction hypothesis for all functions $_{w^+}F'$, defined on
$S^-$ by $_{w^+}F'(s^-)=F(w^+s^-)$, we will show that $_{w^+}F'$
are pluriharmonic. This implies for their Fourier transforms to
satisfy a differential equation with initial data
$f(\lambda,\a,\b)$. Then smoothness of the Fourier transform will
force this function to be zero for $\lambda\notin \bar\Omega \cup
- \bar\Omega$.

Our main work will be to show the smoothness of Fourier
transforms, and will ask for many technicalities.

\bigskip

\noindent {\sl Step 1:  $F'$ is a smooth function of arbitrary order
on $S'$.}
\bpf  First, let $W$ be a right-invariant differential operator on $\nfi
$. We know from the assumptions on $f$ that $Wf$ is well defined, and
bounded.
Therefore, we have
$$
WF'_{s'}(\z ,x)=\int _{\nfi^-}Wf ((\z ,x)s'w^-(s')^{-1})
P'_L(w^-)\ dw^-.
$$
Moreover, partial derivatives of $f$ grow at most
polynomially. The action of $s'$ is linear, hence there are
constants
$C(\a , K)$ and $M(\a )$ such that
$$
|\p _{s'}^{\a }f((\z ,x)s'w^-(s')^{-1})|\leq C(\a ,K)
(1+\tau ( w^-) )^{M(\a )}
$$
for $(\z ,x)s'$ belonging to a compact set $K\subset S'$, with $\tau$
any left--invariant distance as in  (\ref{int}). Now we
select $\eta$ such that $P_{L}$ integrates the right hand side above, to
obtain
$$
\int _{\nfi ^-}|\p _{s'}^{\a }f ((\z,x)s'w^-(s')^{-1})|
P'_L(w^-)\ dw^-<\infty
$$
which allows to  differentiate $F'$ with respect to $s'$.
\epf

\noindent {\sl Step 2:   the function $_{w^+}F'$ satisfies the induction
hypothesis on $S^-$.}
\bpf
We claim first that the assumption (\ref {H2*}), with $S^-$ in place
of $S$, is satisfied for almost
every $w^+$. Indeed, it is sufficient to prove that
\begin{equation}
    \sup _{s'\in S_0'}\| F' _{s'}\| _{L^2(\nfi )}<\infty.
\label{H2'}\end{equation}
 This follows from the fact that,
 for every $s'\in S'_0$, the function $F(\cdot s'\exp t H_{r})$ has a weak limit  in
$L^2(\nfi )$ when $t$ tends to $-\infty$.
Indeed, for $\f \in L^2(\nfi )$,
\begin{align*}
I=&\int _{\nfi }(F(ws'\exp (t_{1}H_{r}))-F(ws'\exp
(t_{2}H_{r})))\f (w)\ dw\\ =&\int _{\nfi } \int _{\nfi }(f(wv_{1})
-f(wv_{2}) P_L(v) \f (w)\ dv\ dw
\end{align*}
with $v_j=s'\exp (t_{j}H_{r})v(s'\exp ( t_{j}H_{r}))^{-1}$ for
$j=1,2$.
 Integrating with respect to $v$ over a compact set $K$ and over its complement
we
get
$$
I\leq \sup _{v\in K}\| f(\cdot v_1)-f(\cdot v_2)\| _{L^2(\nfi )}
\| \f \| _{L^2(\nfi )}+2\| f \| _{L^2(\nfi )}\| \f \| _{L^2(\nfi )}
\int _{K^c}P_{L}(v)\ dv,
$$
which tends to zero when $t_{1},t_{2}\to -\infty$.

\bigskip

We now prove that the functions  $_{w^+}F'$ satisfy the condition
(\ref{equ}), again with $S^-$ in place of $S$. Without loss of
generality, we may assume that $w^+$ is the unit element.
Notice
that the operators
$\sL _j, \De _j$ and $\De ^{\a}_{kj}$
have a perfect sense as left-invariant operators on $S^-$ as far as
indices are smaller than $r$ . Let
$$
(\HH' _j)^T= 2 \De _j + \sum _{k<j}\sum _{\a}\De ^{\a}_{kj}+
    \sum _{j<k<r}\sum _{\a}\De ^{\a}_{jk}.
$$ Again  $(\HH' _j)^T$ may be considered as operators both on $S$
and $S^-$. In the second case $(\HH' _1)^T,...,(\HH' _{r-1})^T$
are $HW$ operators for the tube $V^-+i\W ^-$. We want to prove
that, for $j=1,\cdots, r-1$, we have $$\sL _j F'=(\HH ' _j)^TF'=0.
$$ Since for $i<j<r$, $\sL _j$, $\De^\a_{ij} $, and $\De _j$
commute with $A^+$, we have, for $g'\in S'$,
\begin{align*}
\lim _{t\to -\infty}\sL _jF(g'\exp tH_{r})=&\sL _jF'(g')\\
\lim _{t\to -\infty}\De _jF(g'\exp tH_{r})=&\De _jF'(g')\\
\lim _{t\to -\infty}\De ^\a_{ij} F(g'\exp tH_{r})=&\De ^\a_{ij} F'(g'),
\end{align*}
By hypothesis, $F$ satisfies (\ref{equ}). So we conclude directly
for $\sL _j F'$, $j=1,\cdots, r-1$. For $(\HH ' _j)^TF'$, we
conclude also once we know that
\begin{equation}
    \lim _{t\to -\infty}\De ^{\a }_{jr}F(g'\exp tH_{r})=0
\label{lim2}\end{equation}
Before doing it, we give a last definition. We note $\tilde X^{\a }_{jk}$, $\tilde Y^{\a
}_{jk}$, and $\tilde \sx^{\a }_j$,  $\tilde \sy^{\a }_j$
 the left-invariant vector fields on $N$
which coincide, at the unit element of $N$, with the corresponding elements
of the basis of  ${\cal N}$
that we constructed in subsection 2.3. We define as well
$
\tilde \sL_r=\sum _{\a }(\tilde \sx^{\a }_r)^2+(\tilde \sy^{\a
}_r)^2$.

In the next computation, we identify an element $a$ with a n-uple
$(a_{1}, a_{2}, \cdots a_{r})$, with $a_{j}>0$, in such a way that $a$
is the exponential of $\sum_{j}(\log a_{j})H_{j}$. In particular, an
element $a^+\in A^+$ identifies with a scalar, which we note $a_{r}$
for comprehension. With these notations, the previous limits are
obtained for $a_{r}$ tending to $0$.

Then, it follows from the fact that $ \sL_rF=0$ and a direct
computation  that $$ \x \p _{a_{r} }F(g'a_{r})=\tilde \sL
_rF(g'a_{r}). $$ Moreover,
 $$ \De ^{\a}_{jr}F(g'a_{r})= a_{r}¥
\Big ( \sum _{j<r}\sum_{\a} a_j(\tilde X^{\a }_{jr})^2
+a_j^{-1}(\tilde Y^{\a }_{jr})^2-\frac{1}{\x } \tilde \sL _r \Big
)F(g'a_r)\to 0 $$ when $a_r\to 0$. This finishes the proof of
(\ref{lim2}), as well as the claim of this step. Indeed, for
almost $w^+$, the function $_{w^+}F'$ is pluriharmonic as a
function on $S'$. It follows that $\Delta_{1}F'$ vanishes
identically. This is the main point which will be used later. \epf

\noindent {\sl step 3:   $\hat F'(\lambda, \alpha,\beta, s')$ is a
smooth function of $s'$, for almost every $\lambda $ and every
$\a,\b$.} \bpf As before,  $\hat F'(\lambda, \alpha,\beta, s')$ is
the Fourier transform of the function  $F'_{s'}$, defined on
$N(\Phi)$ by $F'_{s'}(\z, x)=F((\z, x)s')$. We know from (\ref
{H2'}) that it is in $L^2(N(\Phi))$. Moreover, we can write the
Fourier transform of $F'$ in terms of the one of
 the Poisson kernel $P'_L$. Indeed, given $\lambda \in  \L $ and $G_{1},G_{2} \in \sh ^{\lambda }$, we
define the bounded operator $U^{\lambda }_{P'_{s'}}$ by $$
(U^{\lambda }_{P'_{s'}}G_{1},G_{2} )_{\lambda}= \int _{\nfi
^-}P'_{L}(w^-)(U^{\lambda }_{s'w^-(s')^{-1}}G_{1},G_{2})_{\lambda}
dw^-\,,$$ For   $f\in L^1\cap L^2(\nfi )$, it follows directly
from (\ref{P'}) that $$ \hat F'(\lambda, \alpha,\beta,
s')=(U^{\lambda}_f U^{\lambda }_{P'_{s'}} \xi _{\a }^{\lambda
},\xi _{\b }^{\lambda } )_{\lambda} \ \ \text{for}\ a.e. \
\lambda\,. $$
 For
general $f\in L^2(\nfi )$, we use approximation in $L^2(\nfi )$ by
integrable functions and the Plancherel theorem.

So, to prove the claim, it is sufficient to prove the smoothness
$U^{\lambda }_{P'_{s'}}$
  with respect to $s'$. This is given in the following lemma.
\begin{lem}
Assume that $$ \int _{\nfi }\tau(w)^{k+1}P_L(w)\ dw <\infty\,. $$
Then for every $\lambda $ and every $G_{1} ,G_{2}\in \sh
^{\lambda}$ the function $s'\mapsto ( U^{\lambda }_{P'_{s'}}G_{1}
,G_{2})$ is of class ${\cal C}^k$.
\end{lem}
\bpf It follows from (\ref{repr2}) that $(U^{\lambda}_{w^-} G_{1}
,G_{2})$ is a smooth function of $w^-$ with bounded derivatives
(see (6.41) in [OV] and (i), (ii) at the end of (4.1) in [OV]).
Since the action $s'w^-(s')^{-1}$ of $s'$ is linear, $$ |\p ^{\a
}_{s'}(U^{\lambda }_{s'w^-(s')^{-1}}G_{1} ,G_{2})| \leq C(\a
,K)(1+\tau(w^- ))^{|\a |} $$ for  $s'$ belonging to a compact set
$K$. Hence the conclusion follows from the assumption on $P_{L}$.
\epf

This is the end of the step 3.
\epf

\noindent {\sl Step 4:   Conclusion.}

\smallskip
As we said before, we want to write the equation $\Delta_{1}F'=0$ on
the Fourier transform side. We need a preliminary lemma, which will
allow us to do it.

\begin{lem} Let $D$ be an element of the enveloping algebra of
$\sn (\F )^- \oplus (\sn _0)^-\oplus {\cal A}^-$ considered as a
left--invariant operator on $S'$. Then
$$
\sup _{s'\in S_{0}'}\int _{\nfi }|DF'(ws')|^2\ dw <\infty\,.
$$
\label{harn}\end{lem}

\bpf Recall that $_{w^+}F'$ is pluriharmonic on $S^-$. So, by By
the Harnack inequality, we have $$ |DF'(w^+w^-a')|^2 \leq
C(D,B)\int _B |F'(w^+w^-a'g)|^2\ dg, $$ where $B$ is a
neighborhood of identity in $S^-$, and the constant $C(D,B)$ does
not depend on $w^+$.  We use the notation $g=u^-n^-b'$, with
$u^{-}\in N(\Phi)^{-}$, $n^-\in N_{0}^-$ and $b'\in A^-$. Then,
for $s'=ya'\in S'$, we write
\begin{align*}
\int _{\nfi }|DF'&(wya')|^2\ dw =\int _{\nfi }|DF'(y^+w^+w^-y^-
a')|^2\ dw ^-dw^+\\ &\leq C(D,B)\int _B\int _{\nfi
}|F'(y^+w^+w^-y^-a'u^-n^-b')|^2\ dw ^-dw^+du ^-dn^-db'\\ &=
C(D,B)\int _B\int _{\nfi }|F'(y^+w^+w^-y^-a'u^-n^-b' )|^2\  dw
^-dw^+du ^-dn^-db'\\ &= C(D,B)\int _B\int _{\nfi }|F'(wya'n^-b'
)|^2\ dw du ^-dn^-db',
\end{align*}
which is finite. In the above calculation we have used the fact
that the action of $y$ on $\nfi $ is unipotent and we changed
coordinates in $\nfi^-$ in the following way $$
w^-y^-a'u^-(y^-a')^{-1}\to w^-, $$ which preserves the measure
$dw^-$. \epf

\bigskip

We will now prove that for almost every $\lambda$, we have
\begin{equation}
    (-4\pi ^2\< \lambda , \Ad _{s'}X_1 \> ^2 + H_1^2 -H_1)\hat F'(\lambda ,\a ,\b ,s') =0.
\label{equ-four}\end{equation}
To do it, we first approximate $F'$.
 Namely, we take a sequence $\f _n\in {\cal C}_c^{\infty}(\nfi )$
such that $0\leq \f _n \leq 1$, $\f _n=1$ on the ball of radius
$n$, and such that, for every left-invariant differential operator
$D$ on $\nfi $ of positive order,
 $D\f _n\to 0$ uniformly when $n\to \infty$. We put
$$
F'_n((\z ,x)s')=\f _n (\z ,x)F'((\z ,x)s').
$$
A direct calculation shows that
\begin{align*}
\widehat {X_1^mF'_n}(\lambda ,\a ,\b ,s')=& (-2\pi i)^m\< \lambda
,\Ad _{s'}X_1 \> \hat F'_n(\lambda ,\a ,\b ,s'),\\ \widehat{H_1^m
F'_n}(\lambda ,\a ,\b ,s')=& H_1^m\hat F'_n(\lambda ,\a ,\b ,s').
\end{align*}
Then we let $n\to \infty $ and conclude for (\ref {equ-four}) using
Lemma \ref {harn}.
Indeed, Lemma  \ref {harn} implies that
$$
\lim _{n\to \infty} \int _K \int _{\nfi }|DF'_n((\z ,x)s')-DF'((\z
,x)s')|^2\
d\z dx ds'=0
$$
for $D=X_1^2, H_1$ or $H_1^2$ and any compact set $K$ in $S'$.

Now, we prove that (\ref{equ-four}) and the smoothness of $\hat
F'$ forces the Fourier transform of $f$ to vanish outside
$\overline\Omega\cup-\overline\Omega$. Let $J$ be the set of
$\lambda \in V$ such that all the principal minors of $\lambda $
do not vanish (for the definition, see [FK], Proposition VI.3.10).
Since $J$ is dense in $V$, it is sufficient to consider $\lambda
\in J$ such that
$\lambda\notin\overline\Omega\cup-\overline\Omega$ . Then there is
$y_0\in N_0$ such that $\lambda=Ad^\ast_{y_0}\lambda _0$ with
$\lambda _0 = \sum^r_{k=1}b_k\,c_k$, $b_k \neq 0$ for $k=1,...,r$
(see e.g. \cite {DHMP}). Substituting $$ h(\lambda ,s')=\hat F'
(\lambda ,\a ,\b ,y_0^{-1}s') $$ into (\ref{equ-four}), we obtain
$$ \bigl(-4\pi^2\langle \lambda _0,Ad^\ast_{s'}\tilde
X_{1}\rangle^2+H_1^2-H_1\bigr)h(\lambda ,s')=0, $$ or in
coordinates $s'=ya'$, with $a'$ identified with a $r-1$-uple, $$
\bigl(-4\pi^2\langle \lambda _0,Ad^\ast_{y}\tilde
X_{1}\rangle^2+\p ^2_{a_1}\bigr)h(\lambda ,ya')=0, $$ This, and
boundedness of $h$ with respect to $s'$ imply that $$ h(\lambda
,ya') = c(\lambda ,y)\exp {\big [-2\pi |\< \lambda _0,\Ad _yX_1\>
|a_1\big ] }. $$ Letting $a'\to 0$ we get $$ \hat f (\lambda ,\a
,\b )=c(\lambda ,y). $$ Finally, by Lemma 1.27 of \cite{DHMP}
\begin{align}
h(\lambda ,s')=&\hat f (\lambda ,\a ,\b )\exp\big[-2\pi |\<
\lambda _0,\Ad _yX_1\> |a_1 \big ]\notag\\ =& \hat f (\lambda ,\a
,\b )\exp{\big [-2\pi \big |b_1+\frac{1}{2}\sum_{l>1}b_l|y^{1l}|^2
\big |a_1\big ]}\,.\label{end}
\end{align}
Since  for at least one $l>1$ the sign of $b_l$ is different from
the one of  $b_1$, (\ref{end})
 contradicts smoothness of $\hat F'$ with respect to $y$, unless
 $\hat f (\lambda ,\a ,\b )=0$. This concludes for the proof of the lemma.
\epf

\end{document}